\documentclass[12pt]{article}
\newtheorem{prop}{Proposition}[section]
\newtheorem{lem}[prop]{Lemma}

\newtheorem{cor}[prop]{Corollary}
\newtheorem{thm}[prop]{Theorem}
\usepackage{latexsym}
\usepackage{amssymb}
\newcommand{\qbi}[2]{\left[\begin{array}{c}#1\\#2\end{array}\right]}
\newcommand{\pqbi}[2]{\left[\begin{array}{c}#1\\#2\end{array}\right]_{p,q}}

\newcommand{\bu}{{\bf u}}
\newcommand{\bv}{{\bf v}}

\begin{document}

\bibliographystyle{acm}
\title{Set partition statistics and $q$-Fibonacci numbers}

\author{Adam M. Goyt\footnote{This work was partially done while the author was visiting DIMACS.}\\ Department of Mathematics\\ Michigan State University\\ East Lansing, Michigan 48824-1027\\ goytadam@msu.edu\\www.math.msu.edu/$\sim$goytadam \and Bruce E. Sagan\footnotemark[1]\\ Department of Mathematics\\ Michigan State University\\ East Lansing, Michigan 48824-1027\\ sagan@math.msu.edu\\www.math.msu.edu/$\sim$sagan}
\date{August 21, 2006}

\maketitle

\begin{abstract}

We consider the set partition statistics $ls$ and $rb$ introduced by
Wachs and White and investigate their distribution over set partitions
avoiding certain patterns.  In particular, we consider those set
partitions avoiding the pattern $13/2$, $\Pi_n(13/2)$, and those
avoiding both 13/2 and 123, $\Pi_n(13/2,123)$.  We show that the
distribution over $\Pi_n(13/2)$ enumerates certain integer partitions,
and the distribution over $\Pi_n(13/2,123)$ gives $q$-Fibonacci
numbers.  These $q$-Fibonacci numbers are closely related to
$q$-Fibonacci numbers studied by Carlitz and by Cigler.  We provide
combinatorial proofs that these $q$-Fibonacci numbers satisfy
$q$-analogues of many Fibonacci identities.  Finally, we indicate how
$p,q$-Fibonacci numbers arising from the bistatistic $(ls, rb)$ give
rise to $p,q$-analogues of identities.   

\end{abstract}

\section{Introduction and Preliminary Results}

Define the Fibonacci numbers $F_n$ as satisfying the initial
conditions $F_0=1$, $F_1=1$, and the recursion $F_n=F_{n-1}+F_{n-2}$
for $n\ge2$.
This paper focuses on $q$-Fibonacci numbers which arise naturally from
the study of set partition statistics on pattern restricted
partitions.  
Wachs and White~\cite{ww:pqs} introduced the statistics we will use
and showed how that could be used to define $q$-Stirling numbers of
the second kind.
In~\cite{sim:csn}, Simion studied the distribution of these statistics
over non-crossing partitions (those avoiding $13/24$) to get
$q$-analogues of the Catalan numbers.  Further work on this subject
was done by Wachs in~\cite{wac:srg} and by White in~\cite{whi:isp}

The $q$-Fibonacci numbers
studied here are closely related to the $q$-Fibonacci numbers studied
by Carlitz~\cite{car:qfn,car:qfp} 
 and Cigler~\cite{cig:ncq,cig:qfp,cig:aam,cig:qpr}.  In this section
 we provide the necessary 
definitions and background for understanding pattern restricted set
partitions and the two statistics.  Section 2 contains an exploration
of a familiar $q$-analogue of $2^n$ and its relationship to integer
partitions.  In Section 3 we define the $q$-Fibonacci numbers arising
from set partition statistics, introduce the $q$-Fibonacci numbers of
Carlitz and Cigler, and explore their relationships.  The next two
sections are devoted to showing that known Fibonacci identities and
their bijective proofs easily lead to bijective proofs of
$q$-analogues of these identities.  Finally, $p,q$-analogues of
Fibonacci identities are discussed in the last section.  

Let $[n]=\{1,2,\dots,n\}$ and $[k,n]=\{k,k+1,\dots,n\}$.  A {\it
  partition} $\pi$ of $[n]$, written $\pi\vdash[n]$, is a family of
nonempty, disjoint subsets $B_1,B_2,\dots, B_k$ of $[n]$, called {\it
  blocks}, such that $\bigcup_{i=1}^kB_i=[n]$.  If $\pi$ has $k$
blocks then we say that the {\it length} of $\pi$ is $k$, written
$l(\pi)=k$.  We write $\pi=B_1/B_2/\dots/B_k$, omitting set braces and
  commas, and where we always list the
blocks in the {\it standard order} $$\min B_1<\min B_2<\dots<\min
B_k.$$  Suppose $\pi=A_1/A_2/\dots/A_k\vdash[m]$ and
$\sigma=B_1/B_2/\dots /B_l\vdash[n]$.  We say $\pi$ 
{\it is contained in} $\sigma$, written $\pi \subseteq \sigma$, if
  there are $k$  
distinct blocks $B_{i_1},B_{i_2},\dots,B_{i_k}$ of $\sigma$  such that
$A_j\subseteq B_{i_j}$.  For example, if $\sigma=137/25/4/6$ then
$\pi=25/3$ is contained in $\sigma$, but $\pi'=2/5/6$ is not because
the 2 and the 5 would have to be contained in separate blocks of
$\sigma$.   

Given a set of integers $S$ with cardinality $\#S=n$, define the {\it standardization} map \break $St_S:S\rightarrow [n]$ to be the unique order preserving bijection between these sets.  For example if $S=\{3,5,10\}$ then $St_S:S\rightarrow [3]$ and $St_S(3)=1$, $St_S(5)=2$, and $St_S(10)=3$.  We drop the subscript $S$ when this will cause no confusion.  We let $St$ act element-wise on set partitions.  

Let $\pi\vdash[m]$ and $\sigma\vdash[n]$.  We say $\sigma$ contains the {\it pattern} $\pi$ if there is a partition $\pi'$ such that $\pi'\subseteq\sigma$ and $St(\pi')=\pi$, otherwise we say that $\sigma$ {\it avoids} $\pi$.  A copy of $\pi=12/3$ in $\sigma=137/25/4/6$ is $25/6$.  The partition $12/34$ is not contained in $\sigma$ since the only two blocks with more than one element are $\{1,3,7\}$ and $\{2,5\}$, and the latter block can not act as either of the two smallest or two largest elements of 12/34.

Define $$\Pi_n=\{\pi\vdash[n]\},$$ $$\Pi=\bigcup_{n\geq 0} \Pi_n,$$ and for any set of partitions $P\subseteq\Pi$, $$\Pi_n(P)=\{\pi\vdash[n]:\mathrm{\pi \:avoids\: every \:partition\: in\: {\it P}}\}.$$

A {\it layered} partition is a partition of the form $\pi=[1,i]/[i+1,j]/\dots/[k+1,n]$, and a {\it matching} is a partition $B_1/B_2/\dots/B_k$ where $\#B_i\leq 2$ for all $i$.  For example 123/4567/8/9 is a layered partition and 12/34/5/67/8 is a layered matching.  

In~\cite{sag:pas} Sagan characterized $\Pi_n(\pi)$ for
all $\pi\vdash[3]$.  In particular, he showed that $\Pi_n(123)$ is the
set of matchings of $[n]$ and that $\Pi_n(13/2)$ is the set of layered
partitions of $[n]$.  It's not hard to see that
$\#\Pi_n(13/2)=2^{n-1}$.  Goyt~\cite{goy:apt} determined $\Pi_n(P)$
for all sets $P$ of partitions of $[3]$.  He noted that
$\Pi_n(13/2,123)$ is the set of layered matchings of $[n]$, from which
it follows easily that $\#\Pi_n(13/2,123)=F_n$.   

We are interested in the distributions of set partition statistics on $\Pi_n(13/2)$ and $\Pi_n(13/2,123)$ and the resulting $q$-analogues of $2^n$ and $F_n$.  Of the known set partition statistics, only the left smaller, $ls$, and right bigger, $rb$,  statistics of Wachs and White~\cite{ww:pqs} seem to give interesting distributions.  We now describe these statistics. 

Let $\pi=B_1/B_2/\dots/B_k$ be a partition and $b\in B_i$, then we
will say that $(b,B_j)$ is a {\it left smaller pair} of $\pi$ if $j<i$
and $\min B_j<b$.  So, because of standard ordering, for a given block
$B_j$, the elements in left smaller pairs with $B_j$ are exactly those
in $B_{j+1},\dots,B_k$.  We will say that $(b,B_j)$ is a {\it right
  bigger pair} of $\sigma$ if $j>i$ and $\max B_j>b$.  Define
$ls(\pi)$ to be the number of left smaller pairs of $\pi$ and
$rb(\pi)$ to be the number of right bigger pairs of $\pi$.  Wachs and
White proved that $ls$ and $rb$ are equidistributed over
$\Pi_n$.  
turns out that $ls$ and $rb$ are also equidistributed over
$\Pi_n(13/2)$ and $\Pi_n(13/2,123)$.   

\begin{thm}  For any $n$, $$\sum_{\pi\in\Pi_n(13/2)}{q^{ls(\pi)}}=\sum_{\pi\in\Pi_n(13/2)}q^{rb(\pi)},$$
and $$\sum_{\pi\in\Pi_n(13/2,123)}q^{ls(\pi)}=\sum_{\pi\in\Pi_n(13/2,123)}q^{rb(\pi)}.$$
\end{thm}

{\bf Proof:}  Given a set partition
$\pi=B_1/B_2/\dots/B_k\in\Pi_n(13/2)$, let the complement of $\pi$ be
the partition $\pi^c=B_k^c/\dots/B_2^c/B_1^c$, where
$B_i^c=\{n-b+1:b\in B_i\}$.  Notice that taking the complement
reverses the order of the blocks since $\pi$ is layered.  Clearly
complementation is an involution, and so bijective.  To prove the
first equality it suffices to show that it exchanges $ls$ and $rb$.
This follows easily because the block order is reversed and minima are
exchanged with maxima.  Also, complementation does not alter the block
sizes and so restricts to a map on $\Pi_n(13/2,123)$. $\square$   

\section{Distribution over $\Pi_n(13/2)$}

Define $$A_n(q)=\sum_{\pi\in\Pi_n(13/2)}{q^{rb(\pi)}}.$$

It will be useful to think of the $rb$ statistic in the following way.
Consider a block $B_j$.
Let $\pi=B_1/B_2/\dots/B_k$ be a partition.  For each element $b\in
B_i$ with $i<j$ and $b<\max B_j$, we have that $(b,B_j)$ is a right
bigger pair.  The number of right bigger pairs of the form $(b,B_j)$
will be the {\it contribution} of $B_j$ to $rb$.  When restricted to
layered partitions, the contribution of $B_j$ is $$\sum_{i<j}\#B_i=\min
B_j-1.$$  

The generating function $A_n(q)$ is closely related to integer partitions.  A {\it partition} \break $\lambda=(\lambda_1,\lambda_2,\dots,\lambda_k)$ {\it of the integer} $d$ is a weakly decreasing sequence of positive integers such that $\sum_{i=1}^k\lambda_i=d$; the $\lambda_i$ are called {\it parts}.  We let $|\lambda|=\sum_{i=1}^k \lambda_i$.  Denote by $D_{n-1}$ the set of integer partitions with distinct parts of size at most $n-1$.  It is well known that $$\sum_{\lambda\in D_{n-1}}q^{|\lambda|}=\prod_{i=1}^{n-1}(1+q^i).$$  For the rest of this chapter we will refer to a set partition as just a partition and an integer partition by its full name.  

For the following proof, it will be more convenient for us to list the
parts of an integer partition in weakly increasing order.  Let
$\phi:\Pi_n(13/2)\rightarrow D_{n-1}$ be the map defined by
$$\phi(B_1/B_2/\dots/B_k)=(\lambda_1,\lambda_2,\dots,\lambda_{k-1}),$$
where $\lambda_j=\sum_{i=1}^j \#B_i$. 

\begin{thm}  The map $\phi$ is a bijection, and for $\pi\in\Pi_n(13/2)$, $rb(\pi)=|\phi(\pi)|$.  Hence, $$A_n(q)=\prod_{i=1}^{n-1}(1+q^i).$$
\end{thm}

{\bf Proof:}  Given
$\lambda=(\lambda_1,\lambda_2,\dots,\lambda_{k-1})$ consider, for
$1\leq j\leq k$, the differences $d_j=\lambda_j-\lambda_{j-1}$, where
$\lambda_0=0$ and $\lambda_k=n$.  If $\lambda\in D_{n-1}$ then we have
$d_j>0$ for all $j$.  And if $\phi(\pi)=\lambda$ then the $d_j$ give
the block sizes of $\pi$.  So a given $\lambda$ determines a unique
sequence of block sizes, and this sequence determines a unique layered
$\pi$.  Thus, $\phi$ is bijective.  Since $\lambda_j=\sum_{i\le j}\#B_i$
is the contribution of $B_{j+1}$ to $rb$ (and $B_1$ makes no
contribution) we have $rb(\pi)=|\phi(\pi)|$ as desired.  $\square$ 

\section{$q$-Fibonacci Numbers Past and Present}

We turn our focus to the distribution of $rb$ over $\Pi_n(13/2,123)$.  As remarked in the introduction, $$\#\Pi_n(13/2,123)=F_n.$$  The distribution of $rb$ over $\Pi_n(13/2,123)$ gives a nice $q$-analogue of the Fibonacci numbers.  Let $$F_n(q)=\sum_{\pi\in\Pi_n(13/2,123)}q^{rb(\pi)}.$$  

\begin{prop} The generating function $F_n(q)$ satisfies the boundary conditions $F_0(q)=1$, $F_1(q)=1$, and the recursion $$F_n(q)=q^{n-1}F_{n-1}(q)+q^{n-2}F_{n-2}(q).$$\end{prop}

{\bf Proof:}  Let $\pi\in\Pi_n(13/2,123)$.  Since $\pi$ is a matching
it must end in a block of size one or of size two.  If $\pi$ ends in a
singleton then the singleton is $\{n\}$, which contributes $n-1$ to
$rb$, and the remaining elements form a partition in
$\Pi_{n-1}(13/2,123)$.  Similarly the doubleton case is counted by the
second term on the right-hand side of the recursion. $\square$
\medskip

We now introduce the $q$-Fibonacci numbers of Carlitz and Cigler and
explore their relationship to the $q$-Fibonacci numbers as defined
above.  Let $BS_n$ be the set of binary sequences $\beta=b_1\dots b_n$
of length $n$ without consecutive ones.  It is well known that
$\#BS_n=F_{n+1}$.  In~\cite{car:qfn,car:qfp},
Carlitz defined and studied a statistic on $BS_n$ as follows.  Let
$\rho:BS_n\rightarrow \mathbb{N}$ be given by
$$\rho(\beta)=\rho(b_1\dots b_n)=b_1+2b_2+\dots +nb_n,$$ and define
$$F_n^K(q)=\sum_{\beta\in BS_{n-1}}q^{\rho(\beta)}.$$  Carlitz
showed that $F_n^K(q)$ satisfies $F_0^K(q)=1$, $F_1^K(q)=1$, and
$$F_n^K(q)=F_{n-1}^K(q)+q^{n-1}F_{n-2}^K(q).$$   

Cigler~\cite{cig:qfp} defined his $q$-Fibonacci polynomials using
Morse sequences.  A {\it Morse sequence} of length $n$ is a sequence
of dots and dashes, where each dot has length 1 and each dash has
length 2.  For example, $\nu=\bullet\bullet--\bullet-$ is a Morse
sequence of length 9.  Let $MS_n$ be the set of Morse sequences of
length $n$.  Each Morse sequence corresponds to a layered matching
where a dot is replaced by a singleton block and a dash by a
doubleton.  So, $\#MS_n=F_n$.      

Define the weight of a dot to be $x$ and the weight of a dash to be
$yq^{a+1}$ where $a$ is the length of the portion of the sequence
appearing before the dash.  Also, define a weight on the Morse
sequences, $w:MS_n\rightarrow \mathbb{Z}[x,y,q]$, by letting $w(\nu)$
be the product of the weights of its dots and dashes.  For example,
the sequence above has weight
$(x)(x)(yq^3)(yq^5)x(yq^8)=x^3y^3q^{16}$.  Let
$$F_n^C(x,y,q)=\sum_{\nu\in MS_n}w(\nu).$$  Cigler shows that
$F_n^C(x,y,q)$ satisfies $F_0^C(x,y,q)=1$, $F_1^C(x,y,q)=x$, and
$$F_n^C(x,y,q)=xF_{n-1}^C(x,y,q)+yq^{n-1}F_{n-2}(x,y,q).$$  Note that
$F_n^C(1,1,q)=F_n^K(q).$  In fact,
Cigler~\cite{cig:ncq,cig:qfp,cig:aam,cig:qpr} studied more general
$q$-Fibonacci numbers satisfying the above recursion with $yq^{n-1}$
replaced by $t(yq^{n-1})$, where $t$ is an arbitrary nonzero function.
One could apply our method to such $q$-analogues, but we choose $t$ to
be the identity for simplicity.   

\begin{prop}  For all $n\geq0$, $$F_n(q)=q^{{n\choose 2}}F_n^K(1/q).$$
\end{prop}

{\bf Proof:}  It suffices to construct a bijection $\Pi_n(13/2,123)\leftrightarrow BS_{n-1}$ such that if $\pi\leftrightarrow\beta$ then $rb(\pi)={n\choose 2}-\rho(\beta)$.  

Let $\pi\in\Pi_n(13/2,123)$ be mapped to the binary sequence $\beta=b_1 \dots b_{n-1}$ where $b_i=0$ if $i$ and $i+1$ are in separate blocks and $b_i=1$ otherwise.  For example, $1/2/34/56\leftrightarrow 00101$.  We first show that this map is well defined.  Suppose $\pi\mapsto b_1\dots b_{n-1}$, where $b_i=1$ and $b_{i+1}=1$ for some $i$, then $i$, $i+1$, and $i+2$ must be in a block together.  This contradicts the fact that the blocks may only be of size at most 2.  Proving that this map is a bijection is straightforward.  

Now, suppose that $\pi\leftrightarrow\beta$ and $\beta=b_1\dots
b_{n-1}$.  If $\pi=\pi_0=1/2/\dots/n$ so that $b_i=0$ for all $i$,
then
$$
rb(\pi)=\sum_{i=1}^{n-1}i={n\choose2}={n\choose2}-\rho(0\dots0).
$$  
If
$b_i=1$ for some $i$ then $i$ and $i+1$ are in the same block.  In
$\pi_0$, the contribution of the blocks $\{i\}$ and $\{i+1\}$ to $rb$
was $(i-1)+i=2i-1$.  But in $\pi$ the contribution of $\{i,i+1\}$ is
only $i-1$.  Thus, for each $b_i=1$ we reduce $rb(\pi_0)$ by $i$ and
hence,
$$rb(\pi)={n\choose2}-\sum_{i\ :\ b_i=1}i={n\choose2}-\rho(\beta).\:\square$$ 

In order to describe the relationship between $F_n(q)$ and
$F_n^C(x,y,q)$ we will define a weight, $\omega$, on the partitions in
$\Pi_n(13/2,123)$.  Let
$$\omega:\Pi_n(13/2,123)\rightarrow\mathbb{Z}[x,y,q],$$ with
$\omega(\pi)=\omega(B_1/B_2/\dots/B_k)=\prod_{i=1}^k\omega(B_i)$,
where $$\omega(B_j)=\left\{\begin{array}{ll}xq^{\min B_j-1}&
\mathrm{if\:} \#B_j=1,\\ yq^{\min B_j-1}
&\mathrm{if\:}\#B_j=2.\end{array}\right.$$  Now, let
$$F_n(x,y,q)=\sum_{\pi\in\Pi_n(13/2,123)}\omega(\pi).
$$  
Let $s(\pi)$
be the number of singletons of $\pi$, $d(\pi)$ be the number of
doubletons of $\pi$.  It is easy to see directly from the definitions
that
$$F_n(x,y,q)=\sum_{\pi\in\Pi_n(13/2,123)}x^{s(\pi)}y^{d(\pi)}q^{rb(\pi)}.$$
The proof of Proposition 3.1 also shows that 
\begin{equation}
\label{Fn}
F_n(x,y,q)=xq^{n-1}F_{n-1}(x,y,q)+yq^{n-2}F_{n-2}(x,y,q).
\end{equation}

The demonstration of the next result is omitted since it parallels that
of Proposition 3.2, using the bijection
$\Pi_n(13/2,123)\leftrightarrow MS_n$ mentioned above. 

\begin{prop}  For all $n\geq0$, $$F_n(x,y,q)=q^{{n\choose2}}F_{n}^C(x,y,1/q).\:\square$$\end{prop}

\section{$q$-Fibonacci Identities}

We now provide bijective proofs of $q$-analogues of Fibonacci
identities.  Many of the proofs in this paper are simply $q$-analogues
of the tiling scheme proofs of Fibonacci identities given
in~\cite{bq:prc, bcg:tsf}.  
We should also note that Shattuck and Wagner~\cite{sw:pts,sw:nsl} have used various
statitics on domino arrangements to obtain $q$-identities and parity
results for Fibonacci and Lucas numbers.

It is impressive that merely using
the $rb$ statistic on $\Pi_n(13/2,123)$ gives so many identities with
relatively little effort.  We will state our identities for
$F_n(x,y,q)$, but one can translate them in terms of 
$F_n^K(q)$ or $F_n^C(x,y,q)$  using Propositions 3.2 or 3.3, respectively.

\begin{thm}  For all $n\geq0$, $$F_{n+2}(x,y,q)=x^{n+2}q^{{n+2\choose2}}+\sum_{j=0}^nx^jyq^{j+1\choose2}F_{n-j}(xq^{j+2},yq^{j+2},q).$$\end{thm}

{\bf Proof:}  There is exactly one partition in $\Pi_{n+2}(13/2,123)$
with all singleton blocks and the weight of this partition is
$x^{n+2}q^{{n+2\choose2}}$.  The remaining partitions have at least
one doubleton.  Consider all partitions where the first doubleton is
$\{j+1,j+2\}$.  There are exactly $j$ singletons preceding this
doubleton contributing $x^jq^{{j\choose2}}$ to the weight of each
partition.  The doubleton contributes weight $yq^j$.  The remaining
blocks of these partitions form layered matchings of $[j+3,n+2]$.  We
may think of these as being layered matchings of $[n-j]$ where the
contribution of each block to the $rb$ statistic is increased by
$j+2$.  Thus, these contribute weight $F_{n-j}(xq^{j+2},yq^{j+2},q)$.
Hence, the contributed weight of the partitions whose first doubleton
is $\{j+1,j+2\}$, is
$x^jyq^{j+1\choose2}F_{n-j}(xq^{j+2},yq^{j+2},q)$.  Summing from $j=0$
to $n$ completes the proof.  $\square$

\begin{thm}  For all $n\geq0$, $$F_{2n+1}(x,y,q)=\sum_{j=0}^n xy^jq^{j(j+1)}F_{2n-2j}(xq^{2j+1},yq^{2j+1},q),$$ and
$$F_{2n}(x,y,q)=y^nq^{n(n-1)}+\sum_{j=0}^{n-1}xy^jq^{j(j+1)}F_{2n-2j-1}(xq^{2j+1},yq^{2j+1},q).$$\end{thm}

{\bf Proof:}  If $\pi\in\Pi_{2n+1}(13/2,123)$, then $\pi$ must have at least one singleton.  Consider all partitions with first singleton $\{2j+1\}$.  This block must be preceded by $j$ doubletons, which contribute $y^jq^{j(j-1)}$ to the weight.  The singleton $\{2j+1\}$ contributes $xq^{2j}$ to the weight.  The remaining $2n-2j$ elements form a layered matching of $[2j+2,2n+1]$.  As in the previous proof, we may think of these as being elements of $\Pi_{2n-2j}(13/2,123)$ where the contribution of each block to $rb$ is increased by $2j+1$.  This portion of our partition will contribute $F_{2n-2j}(xq^{2j+1},yq^{2j+1},q)$ to the weight.  This proves the first identity.  The proof of the second identity is similar, and hence omitted. $\square$

\begin{thm} For all $n\geq0$ and $m\geq0$, $$F_{m+n}(x,y,q)=F_m(x,y,q)F_n(xq^m,yq^m,q)+yq^{m-1}F_{m-1}(x,y,q)F_{n-1}(xq^{m+1},yq^{m+1},q).$$\end{thm}

{\bf Proof:}  Every $\pi\in\Pi_{m+n}(13/2,123)$ does or does not have $\{m,m+1\}$ as a block.  If $\pi$ has $\{m,m+1\}$ as a block then the blocks prior to this block form a partition in $\Pi_{m-1}(13/2,123)$.  This contributes $F_{m-1}(x,y,q)$ to the weight.  The doubleton $\{m,m+1\}$ has weight $yq^{m-1}$.  The remaining blocks form a partition in $\Pi_{[m+2,m+n]}(13/2,123)$.  The contribution of each block in this partition to the $rb$ statistic is increased by $m+1$, so this portion of the partition contributes $F_{n-1}(xq^{m+1},yq^{m+1},q)$ to the weight.  Thus, the sum of $\omega(\pi)$ over all $\pi$ in $\Pi_{m+n}(13/2,123)$ with doubleton $\{m,m+1\}$ is the second term in the sum above.  

If $\pi\in\Pi_{m+n}(13/2,123)$ does not have $\{m,m+1\}$ as a block, then we can split $\pi$ into a partition of $[m]$ and a partition of $[m+1,m+n]$.  By a similar argument, the sum of $\omega(\pi)$ over all these $\pi$ is $F_m(x,y,q)F_n(xq^m,yq^m,q)$.  $\square$

\medskip

For the proof of the next theorem we will need shifted partitions.  A
partition $\pi\vdash[n]$ {\it shifted} by $k$ positions, denoted
$\pi'$, consists of a block of $k$ blank positions followed by the partition of
$[k+1,k+n]$ obtained by adding $k$ to every element of $\pi$.  For
example, shifting $\pi=134/25$ by $2$ positions gives 
$\pi'=\_ \hspace{3pt}\_ /356/47$.

Let $\Pi_{n,k}(13/2,123)$ be the set of partitions $\pi\in\Pi_n(13/2,123)$ shifted $k$ positions.  Notice that the contribution of each block of $\pi'$ to $rb$ is the contribution of the corresponding block of $\pi$ increased by $k$.    That is $$\sum_{\pi\in\Pi_{n,k}(13/2,123)}\omega(\pi)=F_n(xq^k,yq^k,q).$$

\begin{thm}  For all $m,n\geq1$,

\medskip

\noindent $F_{m+1}(x,y,q)F_{n+1}(xq^m,yq^m,q)
\! = \! xq^mF_{m+n+1}(x,y,q)+
y^2q^{2m-1}F_{m-1}(x,y,q)F_{n-1}(xq^{m+2},yq^{m+2},q).$
\end{thm}  

{\bf Proof:}  The left-hand side is the generating function for all
pairs
$(\pi_1,\pi_2)\in\Pi_{m+1}(13/2,123)\times\Pi_{n+1,m}(13/2,123)$.  The
pair $(\pi_1,\pi_2)$ takes one of two forms.  Either $\pi_1$ ends in a
doubleton and $\pi_2$ begins with a doubleton or not.   

Suppose $\pi_1$ ends in a doubleton and $\pi_2$ begins with a
doubleton.  These doubletons contribute $y^2q^{2m-1}$ to
$\omega(\pi_1)\omega(\pi_2)$.  Dropping these doubletons gives a pair
$(\pi_1',\pi_2')\in\Pi_{m-1}(13/2,123)\times\Pi_{n-1,m+2}(13/2,123)$.
Thus, summing the weights of these pairs gives the second term in the
sum above.   

Suppose $\pi_1$ ends in a singleton, $\pi_2$ begins with a singleton,
or both.  Such a singleton contributes $xq^m$ to
$\omega(\pi_1)\omega(\pi_2)$.  Removing one singleton (in the case of
two singletons, it does not matter which) and concatenating $\pi_1$
and $\pi_2$ gives a partition $\pi\in\Pi_{m+n+1}(13/2,123)$.  Each
$\pi\in \Pi_{m+n+1}(13/2,123)$ has $m+1$ in a block with $m+2$, in a
block with $m$, or in its own block, corresponding to just $\pi_1$
ending in a singleton, just $\pi_2$ beginning with a singleton, or
both.  So every $\pi\in\Pi_{m+n+1}(13/2,123)$ can be constructed as
above.  Thus these contribute weight $xq^mF_{m+n+1}(x,y,q)$.  $\square$ 

\begin{thm}  For all $n\geq0$, $$F_n(x,y,q)F_{n+1}(x,y,q)=\sum_{j=0}^nxy^jq^{\left\lfloor\frac{j^2}{2}\right\rfloor}F_{n-j}(xq^j,yq^j,q)F_{n-j}(xq^{j+1},yq^{j+1},q).$$\end{thm}

{\bf Proof:}  Consider a pair
$$(\pi_1,\pi_2)\in\Pi_n(13/2,123)\times\Pi_{n+1}(13/2,123)
$$ 
with $\pi_1=A_1/A_2/\dots/A_\ell$, and $\pi_2=B_1/B_2/\dots/B_m$.
Search through the blocks in the order $B_1,A_1,$ $B_2,A_2,\dots$ and
find the first singleton block.   Such a block must exists since
either $n$ or $n+1$ is odd.

If the first singleton is some $A_i=\{j\}$ then $B_1,\dots, B_i$ are
all doubletons, and $j$ is odd.  There are $(j-1)/2$ doubletons at the
beginning of $\pi_1$ and $(j+1)/2$ doubletons at the beginning of
$\pi_2$ contributing $y^jq^{(j-1)^2/2}$ to the weight.  The singleton
block $A_i$ has weight $xq^{j-1}$.  The remaining $\ell-i$ blocks of
$\pi_1$ form a layered matching of $[j+1,n]$
providing a contribution of
$F_{n-j}(xq^j,yq^j,q)$.  The remaining $m-i$ blocks of $\pi_2$ are
layered matching of $[j+2,n+1]$ contributing
$F_{n-j}(xq^{j+1},yq^{j+1},q)$.  So, the weight contributed by all
pairs $(\pi_1,\pi_2)$ with $A_i=\{j\}$ as the first singleton is
$$
xy^j q^{\left\lfloor\frac{j^2}{2}\right\rfloor}
F_{n-j}(xq^j,yq^j,q)F_{n-j}(xq^{j+1},yq^{j+1},q).$$ 

If the first singleton is some $B_i=\{j+1\}$ then $j$ is even and, by
similar arguments, the weight contributed by all such pairs is
exactly the same as the one displayed above.
Summing over both even and odd $j$ gives the desired identity. $\square$

\medskip

The identity 
$$
F_n=\sum_{k\geq0}{n-k \choose k}
$$ 
relates the Fibonacci numbers to the binomial coefficients, where
${n\choose k}=0$ if $k>n$.  To state a
$q$-analogue of this identity, we define the {\it $q$-binomial
  coefficients} to be
$$\qbi{n}{k}=\prod_{i=1}^k\frac{q^{n-i+1}-1}{q^i-1},$$  
where, by analogy with binomial coefficients,
$\qbi{n}{k}=0$ if $k>n$.

Carlitz~\cite{car:qfn} derived the following identity using
algebraic and operator methods.  We will provide an alternate proof
using one of the standard combinatorial interpretations of the
$q$-binomial coefficients.  In particular, let $P_{k,l}$ denote the
set of all integer partitions with at most $l$ parts, each of size at
most $k$.  Then \cite[Proposition 1.3.19]{sta:ec1}
\begin{equation}
\label{box}
\qbi{n}{k}=\sum_{\lambda\in P_{k,n-k}} q^{|\lambda|}.
\end{equation}
It will be convenient to represent each $\lambda\in P_{k,l}$ as a path
$p$ in the integer lattice ${\mathbb Z}^2$ from the origin to $(k,l)$, where
each step of $p$ is one unit North ($N$) or one unit East ($E$).  In this case,
the region between $p$ 
and the $y$-axis is just the Ferrers diagram of the corresponding
$\lambda$.  And the area of this region is $|\lambda|$.

\begin{thm}[Carlitz] For all $n\geq0$,
  $$F_n(x,y,q)=
\sum_{k\geq0}x^{n-2k}y^kq^{{k\choose2}+{n-k\choose2}}\qbi{n-k}{k}.$$\end{thm} 

{\bf Proof:}  Let $\Pi_n^k(13/2,123)$ be the set of
$\pi\in\Pi_n(13/2,123)$ with exactly $k$ doubletons.  Thus, we have
$\Pi_n(13/2,123)=\biguplus_{2k\leq n}\Pi_n^k(13/2,123)$, where
$\biguplus$ is the disjoint union.  This implies that
$$F_n(x,y,q)
=\sum_{k\geq0}x^{n-2k}y^k
\left(\sum_{\pi\in\Pi_n^k(13/2,123)}q^{rb(\pi)}\right).
$$   
So it suffices to show that 
\begin{equation}
\label{qrb}
\sum_{\pi\in\Pi_n^k(13/2,123)}q^{rb(\pi)}
=q^{{k\choose2}+{n-k\choose2}}\qbi{n-k}{k}.
\end{equation}

By equation~(\ref{box}), we will be done if we can find a bijection
$\Pi_n^k(13/2,123)\rightarrow P_{k,n-2k}$ such that if
$\pi\leftrightarrow \lambda$ then 
\begin{equation}
\label{rb}
rb(\pi)=|\lambda|+{k\choose2}+{n-k\choose2}.
\end{equation}
Map $\pi=B_1/\ldots/B_{n-k}$ to the lattice path
$p=s_1,\ldots,s_{n-k}$ where $s_i=N$ or $E$ depending on whether $B_i$
is a singleton or doubleton, respectively.  It is easy to see that the
corresponding $\lambda$ is in $P_{k,n-2k}$ and that this is bijective.

As far as the weights, first consider the contribution to $rb(\pi)$ of
those pairs $(b,B_j)$ where $b=\min B_i$ for a doubleton $B_i$. 
If $|B_j|=1$, then $B_i$ and $B_j$ contribute an $E$-step followed later by
an $N$-step in $p$.  Such pairs of steps  are in bijection with squares of the
Ferrers diagram, and thus such $(b,B_j)$ account for the $|\lambda|$ term
of~(\ref{rb}).  If, on the other hand, $|B_j|=2$ then there are
${k\choose2}$ choices for the pair $(b,B_j)$, giving the next term of
our sum.

Finally we need to account for the pairs $(b,B_j)$ where where $b$ is
not the minimum of a doubleton.  But then $b$ could come from any of
the $n-k$ blocks, picking the only element if it is a singleton and
the non-minimum if it is a doubleton.  So there are ${n-k\choose2}$
ways to pick the pair, finishing the proof.  $\square$

\medskip

\begin{thm}  For all $n\geq0$, 
$$F_{2n}(x,y,q)=\sum_{k=0}^n
x^{n-k}y^k   q^{{n+k\choose2}-nk}  \qbi{n}{k}  F_{n-k}(xq^{n+k},yq^{n+k},q).$$
\end{thm}

{\bf Proof:} Let $\Delta_k$ be the set of partitions $\pi\in\Pi_{2n}(13/2,123)$, which begin with a partition of $[n+k]$ having exactly $k$ doubletons.  
Then  $\Pi_{2n}(13/2,123)$ is the disjoint union of the
$\Delta_k$ since the first $n$ blocks of any layered matching of
$[2n]$ must form a partition of the desired type.
Now the same technique used to prove equation~(\ref{qrb}) yields
$$\sum_{\pi\in  \Delta_k}\omega(\pi)=
q^{{n+k\choose2}-nk}x^{n-k}y^k\qbi{n}{k}F_{n-k}(xq^{n+k},yq^{n+k},q).
$$
Summing over all $k$ completes the proof. $\square$

\medskip

We conclude this section by finding a $q$-analogue of the identity
$$
2^n=F_{n+1}+\sum_{k=0}^{n-2}F_k2^{n-2-k}.
$$  
We provide a proof for a $q$-analogue involving only $F_n(q)$, since
we will need to consider blocks with more than 2 elements.     

\begin{thm} For all $n\geq 0$,  
\begin{equation}
\label{prod}
\prod_{i=1}^n(1+q^i)=
F_{n+1}(q)
+\sum_{k=0}^{n-2}q^{k}F_{k}(q)\prod_{i=k+3}^n(1+q^i).
\end{equation} \end{thm} 

{\bf Proof:}  From Theorem 2.1 we have 
$$
\prod_{i=1}^n(1+q^i)=\sum_{\pi\in\Pi_{n+1}(13/2)}q^{rb(\pi)}.
$$  
So we need to show that the right-hand side of equation~(\ref{prod})
also counts $\Pi_{n+1}(13/2)$ with respect to $rb$.

The first term on the right-hand side counts those 
$\pi\in\Pi_{n+1}(13/2)$ that are matchings.  
For any other $\pi$,
suppose the first block of size 3 or larger has
minimum element $k+1$.  The first $k$ elements form a layered matching
of $[k]$, and are hence counted by $F_{k}(q)$.  The block
containing $k+1$ contributes $q^{k}$.  And the remaining blocks contribute
$\prod_{i=k+3}^{n}(1+q^i)$.  $\square$ 

\section{Determinant Identities}

In~\cite{cig:qfp}, Cigler proved a $q$-analogue of the
Euler-Cassini identity, 
$$
F_nF_{n+m-1}-F_{n-1}F_{n+m}=(-1)^nF_{m-1}.
$$
We state his theorem without proof since it will follow from our
results later on in this section.

\begin{thm}[Cigler] For all $n,m\geq1$, the $q$-Fibonacci polynomials $F_n^C(x,y,q)$ satisfy 
$$
F^C_n(x,y,q)F^C_{n+m-1}(x,yq,q)-
F^C_{n-1}(x,yq,q)F^C_{n+m}(x,y,q)
=(-y)^{n}q^{n+1\choose2}F^C_{m-1}(x,yq^{n+1},q).\: \square$$
\end{thm}

Cigler proves this identity twice, once by using determinants and once
by adapting a bijective proof of Zeilberger and
Werman~\cite{wz:bpc}.  We will prove a $q$-analogue of the
Euler-Cassini identity for $F_n(x,y,q)$ by using weighted lattice
paths and their relationship to minors of a Toeplitz-like matrix for the
$q$-Fibonacci sequence.  This is a method that appeared in a paper of 
Lindstr\"{o}m~\cite{lin:vri}, and which was later shown to have broad
application by Gessel and Viennot~\cite{gv:bdp}.   We should note that
Benjamin, Cameron, and Quinn~\cite{bcq:fdc} have recentlt used this
technique to investigate determinants involving ordinary Fibonacci numbers.

Consider the digraph $D=(V,A)$ where the vertices are labeled
$0,1,2,\dots$, and the only arcs are from vertex $n$ to vertex $n+1$
and from vertex $n$ to vertex $n+2$ for all nonnegative integers $n$.
The portion of this digraph consisting of the vertices
$0,1,2,\dots,7$ is pictured below.  All arcs are directed to the
right.   

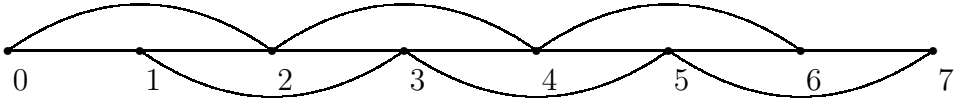
\begin{picture}(400,100)
\put(5,50){\circle*{3}} \put(55,50){\circle*{3}}
\put(105,50){\circle*{3}} \put(155,50){\circle*{3}}
\put(205,50){\circle*{3}} \put(255,50){\circle*{3}}
\put(305,50){\circle*{3}} \put(355,50){\circle*{3}}

\put(7,35){0} \put(57,35){1} \put(107,35){2} \put(157,35){3}
\put(207,35){4} \put(257,35){5} \put(307,35){6} \put(357,35){7}

\put(5,50){\line(1,0){50}} \put(55,50){\line(1,0){50}} \put(105,50){\line(1,0){50}} \put(155,50){\line(1,0){50}} \put(205,50){\line(1,0){50}} \put(255,50){\line(1,0){50}} \put(305,50){\line(1,0){50}}

\qbezier(5,50)(55,85)(105,50) \qbezier(55,50)(105,15)(155,50) \qbezier(105,50)(155,85)(205,50) \qbezier(155,50)(205,15)(255,50) \qbezier(205,50)(255,85)(305,50) \qbezier(255,50)(305,15)(355,50)

\end{picture}

It is easy to see that the number of directed paths from $a$ to $b$ in
$D$ is $F_{b-a}$.  Let the arc from $n$ to $n+1$, written
$\vec{e}_{n,n+1}$, have weight $\omega(\vec{e}_{n,n+1})=xq^n$.  Let
the arc from $n$ to $n+2$ have weight $\omega(\vec{e}_{n,n+2})=yq^n$.
Let $p$ be a directed path from $a$ to $b$, written
$a\stackrel{p}{\rightarrow}b$.  We define the weight of $p$, 
$\omega(p)$, to be the product of the weights of its arcs.  It follows
easily from the definitions that
$$\sum_{p}\omega(p)=F_{b-a}(xq^a,yq^a,q),$$ where the sum is over all
paths $p$ from $a$ to $b$.     

Suppose that $\bu:u_1<u_2<\dots<u_k$ and $\bv:v_1<v_2<\dots<v_k$ are
sequences of  vertices in $D$.  A {\it $k$-tuple}  of
paths from $\bu$ to $\bv$ is
\begin{equation}
\label{P}
P=\{u_1\stackrel{p_1}{\rightarrow}v_{\alpha(1)},u_2\stackrel{p_2}{\rightarrow}v_{\alpha(2)},\dots,u_k\stackrel{p_k}{\rightarrow}v_{\alpha(k)}\}
\end{equation}
where $\alpha\in S_k$, the symmetric group on $k$ elements.  We will let the weight of
such a $k$-tuple be $\omega(P)=\prod_{i=1}^k\omega(p_i)$.  Let
$\mathrm{sgn}(P)=\mathrm{sgn}(\alpha)$, where $\mathrm{sgn}$ denotes
sign.

Now consider the Toeplitz-like matrix  
$$
F=\left[\begin{array}{ccccc} F_0(x,y,q)&F_1(x,y,q)&F_2(x,y,q)&F_3(x,y,q)&\cdots\\
0&F_0(xq,yq,q)&F_1(xq,yq,q)&F_2(xq,yq,q)&\cdots\\
0&0&F_0(xq^2,yq^2,q)&F_1(xq^2,yq^2,q)&\cdots\\
\vdots&\vdots&\vdots&\vdots&\ddots\\
\end{array}\right],
$$
where we label the rows and columns starting with 0. Let $F_{\bu,\bv}$
be the submatrix of $F$ with rows and columns indexed by the sequences
$\bu$ and $\bv$, respectively.  Directly from the definitions we have that   
\begin{equation}
\label{detF}
\det F_{\bu,\bv}=\sum_{P}\mathrm{sgn}(P)\omega(P)
\end{equation}
where the sum is over all $k$-tuples of paths of the form~(\ref{P}).

But we can simplify this sum further.
We will say that two paths are {\it noncrossing} if they do not share a vertex.  

\begin{thm}  
\label{noncross}
Let $\bu:u_1<u_2<\dots<u_k$ and $\bv:v_1<v_2<\dots<v_k$ be
  vertices in $D$.   Then
$$
\det F_{\bu,\bv}=\sum_{P}\mathrm{sgn}(P)\omega(P)
$$
where the sum is over all non-crossing $k$-tuples of paths from $\bu$
to $\bv$.
\end{thm}
{\bf Proof:}  We prove this by giving a weight-preserving,
sign-reversing involution on the $k$-tuples of paths where at least
one pair of paths cross.  Let $k$-tuple $P$ have a crossing pair of
paths.   Let $p_i$ be the path with smallest index
of any path which crosses another path.  Let $w$ be the
first vertex shared by $p_i$ and another path and let $p_j$ be the
path of smallest index $j>i$, that goes through $w$.  Exchange the
portions of $p_i$ and $p_j$ starting at $w$.  This produces a new
$k$-tuple of paths, $Q$, and it is easy to check that this is an
involution.  Since the weight of the $k$-tuple of paths is just the
product of the weights of all arcs appearing in the $k$-tuple, the
weight is preserved.  Finally, the permutations for $P$ and $Q$ differ
by a transposition, so $\mathrm{sgn}(P)=-\mathrm{sgn}(Q)$.  $\square$ 

\medskip

We will now completely characterize the minors of $F$ obtaining, along
the way, a $q$-analogue of the Euler-Cassini identity.  
Given $\bu$ and $\bv$ and a vertex $c$, we call a 
$k$-tuple $P$ of paths from $\bu$ to $\bv$ {\it reducible at $c$} if
no path in $P$ contains both a vertex less than $c$ and a
vertex greater than or equal to $c$.
A family of $k$-tuples is {\it reducible at $c$} if each $k$-tuple is.
If $u_1<c\le v_k$, then it is easy to see that
the sum of the
signed weights for a reducible family can be expressed as a product
over two smaller families.  So it suffices to consider path families
which are not reducible for any $c$.

Now consider the case when $u_i=v_j$ for some $i,j$.  In any $k$-tuple
$P$, this forces $p_i$ to be the path of length $0$ starting and
ending at $c=u_i=v_j$.   The case of all $k$-tuples $P$ 
which are reducible at $c$ has already been covered.  But if $P$ is not
reducible at $c$ then, since paths can connect integers at most two
apart, $P$ contains exactly  one path $p$ beginning before $c$ and
ending after $c$.  Furthermore, this path must contain vertices $c-1$
and $c+1$.  By adding a new endpoint at $c-1$ and a new initial vertex
at $c+1$, one obtains a bijection between all such $P$ and a family of
paths which is reducible at $c$.  Thus, taking into account
 $\omega(\vec{e}_{c-1,c+1})$ and the sign change that
occurs, we can determine the signed sum of the weights of the paths in
this second case using a reducible family.
Thus when $u_i=v_j$ we can compute the determinant using reducible
families, and so
we will assume from now on that $u_i\neq v_j$ for all $i,j$.

The next lemma severely limits the number of minors of $F$ which can
be nonzero.
For a sequence
of vertices $\bu$ and a nonnegative integer $c$, define
$$
\bu(c)=\mbox{number of $u_i< c$}.
$$
\begin{lem}  
\label{ballot}
Suppose the sequences $\bu$ and $\bv$ consist of distinct vertices.
If $\det F_{\bu,\bv}\neq 0$ then we must have
\begin{equation}
\label{bvc}
0\le \bv(c)-\bu(c)\le 2
\end{equation}
for all $c\ge0$.
\end{lem} 
{\bf Proof:}  
We prove both inequalities by contradiction.  Suppose first that
$\bv(c)-\bu(c)<0$.  Then, in a corresponding $k$-tuple $P$, the number of
paths ending before $c$ is greater than the number of paths beginning
in that interval.  Clearly there is no such $k$-tuple and so 
$\det F_{\bu,\bv}=0$ by equation~(\ref{detF}).

On the other hand, suppose  $ \bv(c)-\bu(c)\ge 3$.  Then there must be
at least three paths in $P$ that contain both vertices less than $c$
and vertices greater than or equal to $c$.  Since adjacent vertices on
a path are at most two apart as integers, it is impossible for these
paths to be nonintersecting.  So $\det F_{\bu,\bv}=0$ by
Theorem~\ref{noncross}. $\square$

\medskip

The first inequality in the lemma says that the sequence obtained by
combining $\bu$ and $\bv$ is a {\it ballot sequence}.  But the second
inequality curtails the number of ballot sequences we need to
consider.  Also, if $\bv(c)-\bu(c)=0$ for some $c$ with 
$u_1<c\le v_k$, then any corresponding noncrossing $k$-tuple is
reducible at $c$.  Thus there is only one 
sequence satisfying the lemma which is also irreducible, namely
\begin{equation}
\label{irr}
u_1 < u_2 < v_1 < u_3 < v_2 < u_4 < v_3 < \ldots < u_k < v_{k-1} < v_k.
\end{equation}
So to complete our characterization of the minors of $F$, we need only
consider these sequences.

\begin{thm}
Let $\bu$ and $\bv$ be as in~(\ref{irr}).  Then
\begin{eqnarray*}
\det F_{\bu,\bv}
&=&(-y)^{\sum_{i=1}^{k-1} \left[ v_i-u_{i+1}+1 \right]} \hspace{5pt}
q^{\sum_{i=1}^{k-1} \left[ {v_i\choose2}-{u_{i+1}-1\choose2} \right]}\\
&&\hspace{-30pt} \cdot F_{u_2-u_1-1}(xq^{u_1},yq^{u_1},q)\hspace{3pt}
F_{v_k-v_{k-1}-1}(xq^{v_{k-1}+1},yq^{v_{k-1}+1},q)\hspace{3pt}
\prod_{i=1}^{k-2} F_{u_{i+2}-v_i-2}(xq^{v_i+1},yq^{v_i+1},q).
\end{eqnarray*}
\end{thm}
{\bf Proof:}
Consider a nonintersecting $k$-tuple $P$ counted by 
$\det F_{\bu,\bv}$.  Then $p_1$ 
starts at $u_1$ and must contain the point $u_2-1$ so as not to
intersect $p_2$.  This part of $p_1$ is counted by the factor
$F_{u_2-u_1-1}(xq^{u_1},yq^{u_1},q)$.  Between $u_2$ and $v_1$, the
nonintersecting condition forces $p_1$ to go through exactly the
points having the same parity as $u_2-1$ and $p_2$ to go through the
others.  One of the two paths then terminates at $v_1$ and the other goes from
$v_1-1$ to $v_1+1$.  So the contribution of these steps to the weight
is
\begin{equation}
\label{yq}
yq^{u_2-1}\ y q^{u_2}\ \cdots\ y q^{v_1-1} = 
y^{v_1-u_2+1} q^{{v_1\choose 2}-{u_2-1\choose 2}}.
\end{equation}

Whichever path continues on from $v_1+1$ must then go through $u_3-1$
to avoid intersecting $p_3$, contributing
$F_{u_3-v_1-2}(xq^{v_1+1},yq^{v_1+1},q)$ to the weight.  Next, $p_3$ and
this path alternate vertices between $u_3$ and $v_2$, giving a weight which
is the same as that in equation~(\ref{yq}) but with all the indices
increased by one.

It is clear that this pattern continues, giving the rest of the terms
of the product.  The sign of $P$ is derived in a similar manner, so we
omit the proof.  $\square$

Taking $k=2$ and the  sequences $\bu:0<1$ and $\bv:n<n+m$ in
the previous theorem immediately gives a $q$-analogue of the
Euler-Cassini Identity.  
\begin{cor} For all $n,\geq1$,
$$
F_{n}(x,y,q)F_{n+m-1}(xq,yq,q)-F_{n-1}(xq,yq,q)F_{n+m}(x,y,q)
=(-y)^nq^{{n \choose 2}} F_{m-1}(xq^{n+1},yq^{n+1},q).
$$ 
\end{cor} 
Using this corollary and the identity
$$
F_n(xq^a,yq^a,q)=q^{{n\choose 2}+na}F_n^C(x,y/q^a,1/q)
$$ 
(which is an easy extension of Proposition 3.3) we obtain 
Ciglers' $q$-analogue in Theorem 5.1.

\section{Other Analogues}

The $q$-Fibonacci numbers that are the focus of this paper come from
two statistics, which are equidistributed over the set
$\Pi_n(13/2,123)$.   The next natural question is whether the bistatistic
$(ls,rb)$ also has nice properties when considered on the set
$\Pi_n(13/2,123)$.  The answer is yes.  
Define
$$F_n(x,y,p,q)=\sum_{\pi\in\Pi_n(13/2,123)}x^{s(\pi)}y^{d(\pi)}p^{ls(\pi)}q^{rb(\pi)}.$$
We will also need the $p,q$-{\it binomial coefficient},
$$\pqbi{n}{k}=\prod_{i=1}^k\frac{p^{n-i+1}-q^{n-i+1}}{p^i-q^i}.$$ 
Note that we have  $F_0(x,y,p,q)=1$, $F_1(x,y,p,q)=x$ and for $n\ge2$
$$
F_n(x,y,p,q)= xq^{n-1} F_{n-1}(xp,yp,p,q) + 
yq^{n-2} F_{n-2}(xp^2,yp^2,p,q).
$$

All of the  demonstrations for the formulas in Section 4 translate in a
straightforward manner to the $p,q$ case.
So we will merely list these $p,q$-identities in the following
table and leave the proofs to the reader. 
In this list, we let $F_n(x,y)=F_n(x,y,p,q)$, and
$F_n(xp^a,yp^a)_{1,1}$ be $F_n(xp^a,yp^a,p,q)$ evaluated at $x=y=1$.

In closing, we should note that there are other ways to obtain
$q$-analogues of Fibonacci numbers which could be studied. 
Simion and Schmidt~\cite{ss:rp} discovered  a restricted set
of permutations which is counted by the 
Fibonacci numbers.  There is also a
restricted set of permutations naturally counted by $F_{2n}$, see the
paper of West~\cite{wes:gtf}.
Given the plethora of permutation statistics, some of these sets
should yield interesting $q$-analogues. 

\vspace{50pt}

\noindent {\bf List of $p,q$-Fibonacci Identities}
\hrule
\medskip

\noindent
$F_{n+2}(x,y)=x^{n+2}(pq)^{{n+2\choose2}}+
\displaystyle\sum_{j=0}^n
x^jy \left(\frac{q}{p}\right)^{{j\choose2}} p^{n(j+1)}q^j F_{n-j}(xq^{j+2},yq^{j+2})$
\medskip
\hrule
\medskip

\noindent $F_{2n+1}(x,y)
=\displaystyle\sum_{j=0}^n
xy^j p^{(2n-j)(j+1)-j}q^{j(j+1)}F_{2n-2j}(xq^{2j+1},yq^{2j+1})$

\medskip
\hrule
\medskip
\noindent $F_{2n}(x,y)=y^n(pq)^{n(n-1)}+
\displaystyle\sum_{j=0}^{n-1}
xy^j p^{(2n-j-1)(j+1)-j}q^{j(j+1)}F_{2n-2j-1}(xq^{2j+1},yq^{2j+1})$
\medskip
\hrule
\medskip
\noindent $F_{m+n}(x,y)=F_m(xp^n,yp^n)F_n(xq^m,yq^m)
+\;yp^{n-1}q^{m-1}F_{m-1}(xp^{n+1},yp^{n+1})F_{n-1}(xq^{m+1},yq^{m+1})$

\medskip
\hrule
\medskip
\noindent
$F_{m+1}(xp^n,yp^n)F_{n+1}(xq^m,yq^m)
=xp^n q^m F_{m+n+1}(x,y)$

\medskip

\hspace{.7in}$+\;y^2 p^{2n-1} q^{2m-1} F_{m-1}(xp^{n+2},yp^{n+2})F_{n-1}(xq^{m+2},yq^{m+2}).$

\medskip
\hrule
\medskip

\noindent
$F_n(x,y)F_{n+1}(x,y)=
\displaystyle\sum_{j=0}^n
xy^j p^{n(j+1)-j(j+3)/2}
q^{\left\lfloor\frac{j^2}{2}\right\rfloor}
F_{n-j}(xq^j,yq^j)F_{n-j}(xq^{j+1},yq^{j+1}).$

\medskip
\hrule
\medskip

\noindent
$F_n(x,y)=\displaystyle\sum_{k\geq0} x^{n-2k}y^k
(pq)^{{n\choose2}-k(n-k)}\qbi{n-k}{k}_{p,q}.$

\medskip
\hrule
\medskip

\noindent
$F_{2n}(x,y)=\displaystyle\sum_{k=0}^n
x^{n-k}y^k   (pq)^{{n+k\choose2}-nk}  
\qbi{n}{k}_{p,q}  F_{n-k}(xq^{n+k},yq^{n+k}).$

\medskip
\hrule
\medskip

\noindent
$\displaystyle\prod_{i=1}^n (1+p^{n-i+1}q^i)=
F_{n+1}(x,y)_{1,1}
+\sum_{k=0}^{n-2} q^{k} F_{k}(xp^{n-k+1},yp^{n-k+1})_{1,1}
\prod_{i=k+3}^n (1+p^{n-i+1}q^i).
$

\medskip
\hrule

\vspace{10pt}

{\it Acknowledgement.}  We would like to thank the anonymous referee
for a careful reading of our paper as well as for suggestions that
improved the exposition.
  
\bibliography{ref}
\end{document}